\documentclass[12pt,reqno]{amsart}
\usepackage[a4paper,textwidth=154mm,textheight=236mm,centering]{geometry}
\usepackage{amsmath,amssymb}
\usepackage[hidelinks]{hyperref}
\hypersetup{pdftitle={Hausdorff dimension of reflected Bedford\textendash{}McMullen carpets},pdfauthor={Vyacheslav Koval},pdfkeywords={Bedford-McMullen carpets, Hausdorff dimension, coordinate reflections}}
\newtheorem{theorem}{Theorem}[section]
\newtheorem{lemma}[theorem]{Lemma}
\newtheorem{corollary}[theorem]{Corollary}
\newcommand{\dimH}{\dim_{\mathrm H}}
\newcommand{\one}{\mathbf 1}
\newcommand{\Rows}{\mathcal J}
\title[Dimension of reflected Bedford--McMullen carpets]{Hausdorff dimension of reflected Bedford--McMullen carpets}
\author{Vyacheslav Koval}
\address{Bernoulli Institute, University of Groningen, The Netherlands}
\email{v.v.koval@rug.nl}
\date{}
\begin{document}
\begin{abstract}
We compute the Hausdorff dimension of Bedford--McMullen carpets in which
the copy placed in each chosen cell may be reflected in either coordinate.
Suppose the grid has $n$ columns and $m$ rows, with $2\le m<n$.
The dimension is at least that of the unsigned carpet on the same cells.
The increase is strict
exactly when some row contains both an upright and an upside-down copy
and some row has a different number of chosen cells from its mirror row.
Left--right reflections have no effect. The largest possible increase is
$\log_m 2-\log_n 2$. We give a dimension formula in terms of the numbers
of rectangles in horizontal strips. Two-by-two matrices compute these
counts and give lower and upper bounds at depth $k$ whose difference is
at most $(\log_m 2-\log_n 2)/k$. On each grid a single carpet attains
this bound for every $k$.
\end{abstract}
\maketitle

\section{Main results}

Divide the unit square into $n$ columns and $m$ rows, with $2\le m<n$,
and choose a nonempty set $D$ of cells. In the construction of $K$, the
copy in each chosen cell may be reflected from left to right, from top
to bottom, or both. The same choices are used at every subdivision. Write
$K_0$ for the carpet on the same cells with all orientations preserved,
and put $s_0=\dimH K_0$.

Number the rows from $0$ at the bottom to $m-1$ at the top. Let $t_j$
count the chosen cells in row $j$, let $\tau(j)=m-1-j$ be its mirror
row, and put $\beta=\log m/\log n$. For a row word $w=j_1\cdots j_k$,
let $N(w)$ count the construction rectangles in the corresponding
level-$k$ horizontal strip. Its digits are read from left to right, from
the coarsest subdivision to the finest, in the coordinates of the
original square.

\begin{theorem}[Effect of coordinate reflections]\label{thm:rigidity}
The Hausdorff dimension of $K$ is at least $s_0$. It exceeds $s_0$
if and only if both of the following conditions hold.
\begin{enumerate}
\item Some row contains a copy that preserves the vertical orientation and one that reverses it.
\item Some row $j$ has a different number of chosen cells from its
mirror row, that is, $t_j\ne t_{\tau(j)}$.
\end{enumerate}
For every choice of reflections the dimension is
\begin{equation}\label{eq:main-limit}
 \dimH K=\lim_{k\to\infty}\frac1k\log_m
             \!\left(\sum_{|w|=k}N(w)^\beta\right).
\end{equation}
Changing any of the left--right orientations leaves this value unchanged.
\end{theorem}

The mixed row in the first condition may lie outside the unequal pair
in the second. Either member of that pair may be empty.

Bedford \cite{Bedford} and McMullen \cite{McMullen} independently
computed the dimension of $K_0$ in 1984. Their formula is
\begin{equation}\label{eq:classical}
 s_0=\log_m\!\left(\sum_{j=0}^{m-1}t_j^\beta\right).
\end{equation}
Indeed, in the classical construction $N(j_1\cdots j_k)=t_{j_1}\cdots
t_{j_k}$. The sum in \eqref{eq:main-limit} then equals
$(\sum_j t_j^\beta)^k$, so every depth gives \eqref{eq:classical}.

The first difference in the counts can occur at depth two. Let $a_h$
and $b_h$ be the numbers of upright and upside-down copies in row $h$,
so that $t_h=a_h+b_h$. Within an upright copy, row digit $j$ selects
row $j$ of the pattern. Within an upside-down copy, it selects
row $\tau(j)$. The strip with row digits $h,j$ therefore contains
\begin{equation}\label{eq:two}
 N(hj)=a_h t_j+b_h t_{\tau(j)}
\end{equation}
rectangles. When the row totals agree across every mirror pair, the
two orientations offer the same number of choices. When each occupied row has
just one vertical orientation, the row digits determine
the orientation at each step. These two observations explain the
equality cases in Theorem~\ref{thm:rigidity}.

\subsection*{The two orientation counts}
Put $\Rows=\{0,\ldots,m-1\}$. A cell $e=(i,j)\in D$ in column $i$ and
row $j$ carries signs $(\sigma_e,\eta_e)\in\{+,-\}^2$. With $R_+(u)=u$
and $R_-(u)=1-u$, the map onto $e$ is
\begin{equation}\label{eq:ifs}
 S_e(x,y)=\left(\frac{i+R_{\sigma_e}(x)}{n},\,
                  \frac{j+R_{\eta_e}(y)}{m}\right).
\end{equation}
The carpet $K$ is the unique nonempty compact set satisfying
$K=\bigcup_{e\in D}S_e(K)$. In row $j$ the signs $\eta_e=+$ and $-$
occur $a_j$ and $b_j$ times, respectively. The vertical coordinate
contracts more slowly and is called the weak coordinate. Interchanging
the axes covers $n<m$. For $n=m$, the open cells are disjoint and the
similarity formula gives $\dimH K=\log_m\#D$.

We record counts as row vectors $(x,y)$, in the order upright,
upside down. Appending $j$ on the right of a row word changes its
count vector to $(x,y)M_j$, where
\begin{equation}\label{eq:matrix}
 M_j=\begin{pmatrix} a_j & b_j \\ b_{\tau(j)} & a_{\tau(j)} \end{pmatrix}.
\end{equation}
An upright parent reads row $j$. An upside-down parent reads row
$\tau(j)=m-1-j$, and its $b_{\tau(j)}$ vertically reversed children
are upright in the original square. With $\one=(1,1)^{\mathsf T}$, the
row sums recover the two-step count in \eqref{eq:two}.
\[
 M_j\one=\binom{t_j}{t_{\tau(j)}},\qquad
 (1,0)M_hM_j\one=a_h t_j+b_h t_{\tau(j)}.
\]
For $w=j_1\cdots j_k$ put $M_w=M_{j_1}\cdots M_{j_k}$. The original
square is upright, so $N(w)=(1,0)M_w\one$.

The mirror word $\tau(w)=\tau(j_1)\cdots\tau(j_k)$ has its digits in
the same order as $w$. Let $P=\left(\begin{smallmatrix}0&1\\1&0\end{smallmatrix}\right)$
exchange the two orientations. Equation~\eqref{eq:matrix} gives
$M_{\tau(j)}=P M_j P$. Since $P^2$ is the identity and $P\one=\one$,
we get $M_{\tau(w)}=P M_w P$ and $(0,1)M_w\one=N(\tau(w))$.
The sums for an upright start and an equally weighted start are
\begin{equation}\label{eq:sums}
 Z_k=\sum_{w\in\Rows^k}N(w)^\beta,\qquad
 U_k=\sum_{w\in\Rows^k}
       \left(\frac{N(w)+N(\tau(w))}2\right)^\beta .
\end{equation}

\begin{theorem}[Dimension formula]\label{thm:pressure}
For every carpet \eqref{eq:ifs},
\begin{equation}\label{eq:pressure}
 \dimH K=\sup_{k\ge1}\frac{\log_m Z_k}{k}
        =\inf_{k\ge1}\frac{\log_m U_k}{k}.
\end{equation}
Both normalized logarithms converge to $\dimH K$ as $k\to\infty$. At every depth $k$,
\begin{equation}\label{eq:error}
 0\le \frac{\log_m U_k-\log_m Z_k}{k}
     \le\frac{(1-\beta)\log_m 2}{k}.
\end{equation}
For each $m<n$ there is a carpet that attains the last bound at every depth.
\end{theorem}

At depth one, Theorem~\ref{thm:pressure} gives the bounds
\begin{equation}\label{eq:extrema}
 s_0\le\dimH K\le s_*,
 \qquad s_*:=\log_m\!\left(\sum_{j=0}^{m-1}
                 \left(\frac{t_j+t_{\tau(j)}}2\right)^\beta\right).
\end{equation}
Section~\ref{sec:proof} proves the dimension formula by comparison with
ordinary carpets on finer grids. Section~\ref{sec:rigid} proves the
criterion for a strict increase and determines when the upper bound $s_*$
is attained. Section~\ref{sec:algorithm} gives bounds that can be
computed without enumerating all $m^k$ row words.

\subsection*{Relation to earlier work}

Kenyon and Peres \cite{KP} computed the Hausdorff dimension of planar
sofic affine-invariant sets. Their theorem applies here after keeping
track of both coordinate orientations on a graph with at most four
vertices. The graph represents a finite union of isometric copies of $K$.
Takahashi \cite{Takahashi} studied dimension spectra in this setting,
Barral and Feng \cite{BF} developed the corresponding weighted pressure
in greater generality, and Alibabaei \cite{Alibabaei} obtained exact
evaluations for certain sofic sets. For the present carpets, the two
vertical states are related by the single involution $\tau$.
We use this relation to obtain the sharp finite-depth estimate
\eqref{eq:error} and the equality cases in \eqref{eq:extrema}.

Signed maps occur naturally in graphs of continuous self-affine functions.
Urba\'nski \cite{Urbanski} computed dimensions in this class, and Allaart
and Kawamura \cite{AK} explain its relation to finite-pattern carpets.
Fraser \cite{Fraser} treated reflections in his work on the box and
packing dimensions of box-like sets. B\'ar\'any and David \cite{BD}
give Hausdorff dimension formulas for weakly separated carpets and discuss
simultaneous vertical reversals within rows. Here a single row may contain
both signs, and distinct vertical words may even give the same affine map.
Section~\ref{sec:branches} characterizes these coincidences and their
effect on the dimension. Section~\ref{sec:box} compares the Hausdorff
and box dimensions.

\section{The classical row sum}\label{sec:classical}

A row address of length $k$ locates the vertical coordinate to within
$m^{-k}$. Only $\lfloor\beta k\rfloor$ column digits are needed for
comparable horizontal precision. McMullen \cite{McMullen} chose a
Bernoulli measure adapted to these approximate squares. Peres gave a
shorter proof of \eqref{eq:classical} by applying the Rogers--Taylor
density theorem to this measure \cite[Sect.~2, pp.~515--517]{Peres}.
We give the approximate-square argument here because the proof for
reflected carpets uses the classical formula on successively finer grids.

Let $Z=\sum_j t_j^\beta$. On the address space $D^{\mathbb N}$ choose cells
independently, with row probabilities $q_j=t_j^\beta/Z$ and equal
probabilities among the cells of each occupied row. Call the resulting
measure $\mu$. For an address $\omega$ with row digits $j_1,j_2,\ldots$,
write $S_k=\sum_{r=1}^k\log t_{j_r}$ and $S_0=0$. An approximate square
$Q_k(\omega)$ specifies the first $k$ row digits and the first
$\ell=\lfloor\beta k\rfloor$ column digits. Its image lies in a grid
rectangle of height $m^{-k}$ and width $n^{-\ell}\in[m^{-k},n m^{-k})$.
Specifying a column in row $j_r$ costs a factor $1/t_{j_r}$. Hence
\begin{equation}\label{eq:classmass}
 \mu(Q_k(\omega))=\prod_{r=1}^k q_{j_r}
                     \prod_{r=1}^{\ell}t_{j_r}^{-1}
                   =Z^{-k}\exp(\beta S_k-S_\ell).
\end{equation}

For $\mu$-almost every address, the law of large numbers gives convergence
of $S_k/k$. Consequently $(\beta S_k-S_\ell)/k\to0$.
Fix $0<s<s_0$. On a set of positive measure the bound
$\mu(Q_k(\omega))\le m^{-ks}$ holds for every $k$ above a common
threshold. Restrict $\mu$ to this set and project it to the carpet.
A ball of radius $m^{-k}$ meets a bounded number of level-$k$ grid
rectangles. Its preimage is contained in the union of the corresponding
symbolic atoms, each of restricted mass at most $m^{-ks}$. Thus the projected measure satisfies the mass distribution bound at
exponent $s$. This proves $\dimH K_0\ge s_0$, with the case $s_0=0$
immediate.

For the upper bound, fix any address. Choose $k$ tending to infinity
along a sequence on which $S_k/k$ approaches its limsup. The averages
$S_\ell/\ell$ have at most the same limsup, so
$\limsup_k (\beta S_k-S_\ell)/k\ge0$.
By \eqref{eq:classmass}, for every $s>s_0$ this address belongs to
arbitrarily fine approximate squares of mass at least $m^{-ks}$.
Starting at any prescribed depth $k_0$, stop each address the first time
this happens. The stopped squares partition address space because the
partitions are nested. Their images cover $K_0$, have diameters at most
$\sqrt{1+n^2}\,m^{-k_0}$, and the sum of their diameters to the power
$s$ is at most $(1+n^2)^{s/2}\sum_Q\mu(Q)=(1+n^2)^{s/2}$.
Letting $k_0\to\infty$ proves $\dimH K_0\le s_0$.

\section{Counting oriented rectangles}\label{sec:proof}

A level-$k$ grid cell is specified by its two geometric digit words, regardless of its orientation. Different words in $D^k$ select different cells, since their first differing letters select different subcells of a common parent. Consequently incidences along cell boundaries do not alter any of the counts. From either initial horizontal orientation the number of choices for a given weak transition is the corresponding entry of \eqref{eq:matrix}.

We first compare $K$ with two ordinary carpets on finer grids.
For the outer approximation we allow the orientation to change freely
between blocks. For the inner approximation we keep blocks of a single
orientation and concatenate them in pairs. There are only four possible
sign pairs, so the loss in the logarithmic row sum is bounded independently
of the length of the block.

\begin{lemma}\label{lem:geometric}
For every $k\ge1$,
\[
 \frac{\log_m Z_k-\log_m4}{k}\le \dimH K
     \le \frac{\log_m Z_k+\log_m4}{k}.
\]
\end{lemma}
\begin{proof}
For the upper bound, allow each $k$-block to start in any of the four coordinate orientations, independently of the preceding block. Keep all geometric grid cells that can occur, and use the positive map onto each of them. The resulting unsigned $(n^k,m^k)$-carpet contains $K$, since every infinite address of $K$ uses only the permitted geometric blocks. This argument also covers points with two grid expansions, since it only requires one address. At most $2N(w)+2N(\tau(w))$ cells lie over $w$. The unsigned formula and $(x+y)^\beta\le x^\beta+y^\beta$ bound the row sum of this carpet by $2^{1+\beta}Z_k\le4Z_k$.

For the lower bound let $N_g(w)$ count the words over $w$ with total sign pair $g$. Since $N(w)^\beta\le\sum_g N_g(w)^\beta$, summing over $w$ shows that some fixed $g$ satisfies
\[
   \sum_w N_g(w)^\beta\ge Z_k/4 .
\]
Take all ordered pairs of length-$k$ words with this sign pair. Every resulting $2k$-block has both signs positive since $g^2$ is the identity. Concatenating such blocks gives an unsigned subcarpet of $K$. For first and second row words $u,v$ read in their own positive frames, the geometric row of the pair is $uv$ or $u\tau(v)$, according to the vertical sign of $g$. The map from $(u,v)$ to that row is injective, and there are $N_g(u)N_g(v)$ rectangles over it. Its row sum is therefore $(\sum_w N_g(w)^\beta)^2$. Applying the unsigned formula at depth $2k$ proves the bound.
\end{proof}

\begin{proof}[Proof of Theorem~\ref{thm:pressure}]
It is useful to allow a mixture of starting orientations. Give the positive state weight $p$ and the negative state weight $1-p$, and put $v(p)=(p,1-p)$. The sum
\begin{equation}\label{eq:F}
 F_k(p)=\sum_{w\in\Rows^k}\big(v(p)M_w\one\big)^\beta
       =\sum_w\big(pN(w)+(1-p)N(\tau(w))\big)^\beta
\end{equation}
is concave in $p$ and symmetric about $1/2$. Its minimum is $Z_k=F_k(1)$ and its maximum is $U_k=F_k(1/2)$. After reading a prefix $u$, the two entries of $v(p)M_u$ count its descendants in each orientation with these initial weights. If their sum $c=v(p)M_u\one$ is positive, let $q$ be the first entry divided by $c$. The relative weights for the next block are then $(q,1-q)$. The contribution of its extensions by $\ell$ digits is $c^\beta F_\ell(q)$. Summing over prefixes gives
\begin{equation}\label{eq:concatenate}
 Z_\ell F_k(p)\le F_{k+\ell}(p)\le U_\ell F_k(p).
\end{equation}
Thus $Z_k$ is supermultiplicative and $U_k$ is submultiplicative. Subadditivity of $x^\beta$ gives $U_k\le 2^{1-\beta}Z_k$. Lemma~\ref{lem:geometric} identifies the common logarithmic growth rate with $\dimH K$. Fekete's lemma gives \eqref{eq:pressure}, and the ratio bound gives \eqref{eq:error}.

For sharpness choose two cells in row $0$ and no others. Place an upright
copy in one cell and an upside-down copy in the other. From an upright parent both children are in geometric
row $0$. From an upside-down parent both are in row $m-1$. In either case
one child is upright and one is upside down. The first digit must therefore
be $0$, and the count vector after that digit is $(1,1)$. Appending either
$0$ or $m-1$ keeps this vector equal to $(1,1)$, since exactly one of
the two parents contributes its two children. No other digit is possible.
Thus there are $2^{k-1}$ occupied row words and each strip contains two
distinct rectangles, one in each orientation. The mirror words all begin
with $m-1$ and are disjoint from the occupied words. It follows that
\begin{equation}\label{eq:sharp-counts}
 Z_k=2^{k-1}2^\beta,\qquad U_k=2^k .
\end{equation}
The difference in \eqref{eq:error} is exactly $(1-\beta)\log_m2/k$
for every $k$.
\end{proof}

\section{When reflections increase the dimension}\label{sec:rigid}

The second row digit in \eqref{eq:two} can be read in either orientation.
The two resulting vectors of counts, $(t_j)_j$ and $(t_{\tau(j)})_j$,
are permutations of one another. A mixed row averages these vectors in
proportions given by its two signs. Since $x^\beta$ is strictly concave,
averaging increases the sum of the powers unless the two vectors agree.
The lower bound in Theorem~\ref{thm:pressure} carries this two-step gain
to the dimension of the carpet.

\begin{proof}[Proof of Theorem~\ref{thm:rigidity}]
The limit formula follows from Theorem~\ref{thm:pressure}. The matrices
in \eqref{eq:matrix} use only vertical orientations, so changing a
left--right orientation leaves every $N(w)$ unchanged. It remains to
decide when the lower bound in \eqref{eq:extrema} is strict.
For one digit the sum in \eqref{eq:F} becomes
\[
 H(p):=F_1(p)=\sum_{j=0}^{m-1}\bigl(p\,t_j+(1-p)t_{\tau(j)}\bigr)^\beta.
\]
If all mirror pairs have equal totals, $M_j\one=t_j\one$. It follows that $N(w)=\prod_i t_{w_i}$ for every $w$, so $\dimH K=s_0=s_*$. Otherwise $H$ is strictly concave because at least one summand is strictly concave. Hence
\[
 H(0)=H(1)<H(p)\quad(0<p<1),\qquad
 H(p)<H(1/2)\quad(p\ne1/2).
\]

From positive orientation, row $j$ has mass $t_j$ and gives the next orientation proportion $a_j/t_j$ when $t_j>0$. Thus
\[
 Z_2=\sum_{t_j>0}t_j^\beta H(a_j/t_j).
\]
Any mixed row gives $Z_2>Z_1^2$, hence $\dimH K>s_0$ by Theorem~\ref{thm:pressure}. With no mixed row, the orientation is determined by each geometric prefix. Every positive-mass prefix starting from $p=1$ has next proportion $0$ or $1$, so $Z_{k+1}=Z_1 Z_k$. This proves the criterion.
\end{proof}

\begin{theorem}[Upper bound and largest increase]\label{thm:extrema}
Write $\delta_j=a_j-b_j$. The upper bound $s_*$ in \eqref{eq:extrema}
is attained precisely when $t_j=t_{\tau(j)}$ for every $j$ or
$\delta_j=\delta_{\tau(j)}$ for every $j$. For prescribed row totals
there is a choice of signs attaining $s_*$ if and only if $t_j$ and
$t_{\tau(j)}$ have the same parity for every $j$. On each grid the largest
possible increase over all choices of cells and reflections is
\begin{equation}\label{eq:gain}
 \max_{D,\,\text{signs}}\bigl(\dimH K-\dimH K_0\bigr)
       =\log_m2-\log_n2 .
\end{equation}
\end{theorem}

\begin{proof}
Equal totals in all mirror pairs give $s_0=s_*$. Suppose now that some pair has unequal totals. Then $H=F_1$
has its unique maximum at $1/2$, as in the preceding proof.

Starting with equal weights on the two orientations gives mass
$s_j=(t_j+t_{\tau(j)})/2$ to digit $j$. If $s_j>0$, the new proportion
of upright rectangles is the first entry of $(1/2,1/2)M_j$ divided by
$s_j$. Thus
\[
 q_j=\frac{a_j+b_{\tau(j)}}{t_j+t_{\tau(j)}},\qquad
 q_j-\frac12=\frac{\delta_j-\delta_{\tau(j)}}{2(t_j+t_{\tau(j)})}.
\]
If some such $q_j$ differs from $1/2$, then
$U_2=\sum_{s_j>0}s_j^\beta H(q_j)<U_1^2$ and
$\dimH K\le\frac12\log_m U_2<s_*$. If every $q_j=1/2$ on the
occupied pairs, the proportion stays at $1/2$ after every possible
word, so $U_k=U_1^k$ and $\dimH K=s_*$. For an empty pair both
imbalances are zero. This proves the first assertion.

When the totals are mirror-symmetric, every sign assignment attains
$s_*$ and the parity condition holds automatically. Otherwise equality
requires equal imbalances in every mirror pair. Since
$\delta_j=t_j-2b_j$, this requires equal parities. Conversely, for
totals of matching parity choose $\delta_j=0$ for even $t_j$ and
$\delta_j=1$ for odd $t_j$. Taking $a_j=(t_j+\delta_j)/2$ and
$b_j=(t_j-\delta_j)/2$ realizes these imbalances in each row.

Finally, $U_1\le 2^{1-\beta}Z_1$ bounds the increase by
$(1-\beta)\log_m2=\log_m2-\log_n2$. In the two-cell example of
\eqref{eq:sharp-counts} the reflected dimension is $\log_m2$ and the
unsigned dimension is $\log_n2$.
\end{proof}

For prescribed row totals, taking all copies upright always attains
$s_0$. Attaining the maximum in \eqref{eq:gain} requires equality in
both $\dimH K\le s_*$ and $U_1\le 2^{1-\beta}Z_1$.
Strict subadditivity of $x^\beta$ shows that the latter equality holds
exactly when $t_jt_{\tau(j)}=0$ for every $j$. By the first part of the
theorem the maximum is therefore attained exactly when the occupied
rows are disjoint from their mirrors and each has equally many copies
of either vertical sign.

\subsection*{An exact family between the two bounds}
The gap between the endpoints in \eqref{eq:extrema} can already be seen
with just two occupied rows. Choose a row $r$ different from its mirror
$\tau(r)$. Put $t$ cells in $r$ and $u$ cells in $\tau(r)$, with
$1\le t,u\le n$ and $t\ne u$. Let $p\in[0,1]$ satisfy
$pt,pu\in\mathbb Z$. Give a proportion $p$ of the cells in $r$ and
a proportion $1-p$ in $\tau(r)$ positive vertical sign. The two row matrices factor as
\[
 M_r=\binom{t}{u}(p,1-p),\qquad
 M_{\tau(r)}=\binom{u}{t}(1-p,p).
\]
Reading a digit now determines the relative orientation counts, whatever
the preceding digits were. Consequently the strip count is multiplied by $A=pt+(1-p)u$ when the digit repeats
and by $B=(1-p)t+pu$ when it changes. Summing over the two possibilities at
each step gives
\[
 Z_k=(t^\beta+u^\beta)(A^\beta+B^\beta)^{k-1},
 \qquad
 \dimH K=\log_m\!\left(A^\beta+B^\beta\right).
\]
At $p=0$ or $1$ this is the classical row sum. At $p=1/2$, when
admissible, it is the mirror-averaged upper bound. For $0<p<1/2$ the
difference $|A-B|=(1-2p)|t-u|$ decreases as $p$ increases while $A+B$
stays fixed. Strict concavity of $x^\beta$ gives strict increase along
the admissible values, and replacing $p$ by $1-p$ exchanges $A$ and $B$.
Thus this family gives an exact value even when the two mirror rows both
contain both signs.

\subsection*{Graphs at the upper endpoint}
For another equality case, consider the graph of a continuous self-affine function from $[0,1]$ to $[0,1]$ joining $(0,0)$ to $(1,1)$. Suppose there is one piece per column and all horizontal orientations are positive \cite{Urbanski,AK}. Write the first-level endpoint heights in units of $1/m$. Moving across the columns gives a nearest-neighbour walk from $0$ to $m$. For each $j$ its net number of upcrossings of $[j,j+1]$ is one, by telescoping the indicator of being above $j$. Thus $a_j-b_j=1$ throughout, and Theorem~\ref{thm:extrema} gives $\dimH K=s_*$.

\section{Computing the dimension by interpolation}\label{sec:algorithm}

The sums $Z_k$ and $U_k$ in Theorem~\ref{thm:pressure} run over all
$m^k$ row words. There is another way to bound their growth rate.
After a row word has been read, its extensions depend only on the total
number of rectangles and the proportion $p$ that are upright. We can
therefore approximate the sum over extensions by interpolating in the
single variable $p$. The method also applies when some row matrices are
singular. For related transfer-operator methods see Morris \cite{Morris}
and Falk and Nussbaum \cite{FN}.

With $v(p)=(p,1-p)$ as in \eqref{eq:F}, put
\[
 \begin{aligned}
 c_j(p)&=v(p)M_j\one=p\,t_j+(1-p)t_{\tau(j)},\\
 T_j(p)&=\frac{(v(p)M_j)_1}{c_j(p)}
        =\frac{p\,a_j+(1-p)b_{\tau(j)}}{p\,t_j+(1-p)t_{\tau(j)}}
          \qquad(c_j(p)>0).
 \end{aligned}
\]
Here $c_j(p)$ is the factor by which the count changes on reading $j$,
and $T_j(p)$ is the new upright proportion. In particular,
$T_j(1/2)=q_j$ whenever $s_j>0$. Thus the functions $F_k$
in \eqref{eq:F} satisfy $F_k=\mathcal L^k1$, where
\[
 (\mathcal L f)(p)=\sum_{j:c_j(p)>0}c_j(p)^\beta f(T_j(p)).
\]
Each term extends continuously by zero where $c_j(p)=0$, so $\mathcal L$
maps continuous functions on $[0,1]$ to continuous functions. The
symmetry and concavity of $F_k$ give $Z_k\le F_k(p)\le U_k$.

Choose nodes $0=p_0<\cdots<p_N=1$. Let $E$ denote restriction to
these nodes and $I$ piecewise linear interpolation of nodal values.
The nonnegative matrix $B=E\mathcal L I$ has order $N+1$. To form its $i$th
row, distribute the weight $c_j(p_i)^\beta$ between the two nodes
bracketing $T_j(p_i)$ in their linear interpolation proportions, for
each $j$ with $c_j(p_i)>0$. Each row has at most $2m$ nonzero entries.
Write $\mathcal P$ for the set of nodes and put
\[
 h(p)=p^\beta+(1-p)^\beta,\qquad
 \gamma(\mathcal P)=\max_{0\le p\le1}\{h(p)-(IEh)(p)\}.
\]

Since $1\le h\le 2^{1-\beta}$, linear interpolation gives
$0\le\gamma(\mathcal P)\le 2^{1-\beta}-1<1$.

\begin{theorem}\label{thm:grid}
For every grid $\mathcal P$,
\[
 \frac{\log\rho(B)}{\log m}\ \le\ \dimH K\ \le\
 \frac{\log\rho(B)-\log(1-\gamma(\mathcal P))}{\log m}.
\]
For each integer $L\ge1$ there is a grid with $2L+1$ nodes for which
$\gamma(\mathcal P)\le C_\beta L^{-2}$, with $C_\beta$ depending
only on $\beta$.
\end{theorem}

\begin{proof}
Each summand in $F_k$ has the form
$f(p)=(pa+(1-p)b)^\beta$, with $a,b\ge0$. When $a>b$,
the function $(a-b)^\beta p^\beta-f(p)$ is concave. Indeed,
$pa+(1-p)b\ge p(a-b)$, and comparison of second derivatives on
$(0,1)$ gives the assertion. Thus
\[
 0\le f-IEf\le (a-b)^\beta\big(p^\beta-IE(p^\beta)\big).
\]
For $b>a$ the same estimate holds with $1-p$ in place of $p$ and
$b-a$ in place of $a-b$. In summing over words we have
$\sum_w a^\beta=\sum_w b^\beta=Z_k$. Hence
\[
 0\le F_k-IE F_k\le Z_k(h-IEh)\le \gamma(\mathcal P) F_k.
\]
Applying $E\mathcal L$ and using positivity gives
$(1-\gamma(\mathcal P)) E F_{k+1}\le B E F_k\le E F_{k+1}$.
Induction bounds $B^k\mathbf 1_{N+1}$ between
$(1-\gamma(\mathcal P))^k E F_k$ and $E F_k$, where $\mathbf 1_{N+1}$ is the
all-ones nodal vector. Since the $k$th roots of $Z_k$ and $U_k$
both tend to $m^{\dimH K}$, the spectral radius formula in the
maximum row-sum norm gives
\[
 (1-\gamma(\mathcal P))m^{\dimH K}\le\rho(B)\le m^{\dimH K}.
\]
Taking logarithms yields the stated interval of width
$-\log(1-\gamma(\mathcal P))/\log m$.

For the rate, put $s=2/\beta$, take $x_i=\tfrac12(i/L)^s$ for
$0\le i\le L$, and include the reflected nodes $1-x_i$. The midpoint
occurs in both halves, giving $2L+1$ distinct nodes. On $[0,x_1]$ the
oscillation of $h$ is at most $x_1^\beta=2^{-\beta}L^{-2}$.
For $1\le i<L$ the second derivative of $h$ on $[x_i,x_{i+1}]$ is
bounded in absolute value by $C_\beta x_i^{\beta-2}$. Since
$(i+1)^s-i^s\le C_s i^{s-1}$, the linear interpolation error there
is at most
\[
 C_\beta (x_{i+1}-x_i)^2x_i^{\beta-2}
 \le C'_\beta L^{-s\beta} i^{s\beta-2}= C'_\beta L^{-2}.
\]
Symmetry gives the same bound on the right half.
\end{proof}

On each grid interval the maximum defining $\gamma(\mathcal P)$ is found by
equating $h'$ to the slope of the secant. If only bounds on $\rho(B)$
are required, every positive nodal vector $x$ gives
$\min_i(Bx)_i/x_i\le\rho(B)\le\max_i(Bx)_i/x_i$.

\section{When can the two signs be counted separately?}\label{sec:branches}

Another possible extension of the classical formula is to count the
upright and upside-down cells in a row as two different choices of
vertical branch. Two different branch
words can describe exactly the same map of the vertical interval.
If such an equality occurs at any depth, there are already two words
of length two that define the same map.

Write $g_{j,+}(y)=(j+y)/m$ and $g_{j,-}(y)=(j+1-y)/m$.
Keep the branches with positive multiplicity, and denote their alphabet
by $\mathcal B$. Thus the weights of $(j,+)$ and $(j,-)$ are $a_j$
and $b_j$, respectively. Put
\[
 U=\sum_{j=0}^{m-1}(a_j^\beta+b_j^\beta).
\]
An \emph{exact overlap} is an equality of affine maps represented by
distinct words over $\mathcal B$. Such words necessarily have the
same length.

The open set condition asks for a nonempty bounded open set $O\subset\mathbb R$ such that the images $g(O)$ are pairwise disjoint subsets of $O$.
\begin{theorem}\label{thm:branches}
One has $\dimH K\le \log U/\log m$. The following are equivalent.
\begin{enumerate}
\item Equality holds.
\item The signed weak branches satisfy the open set condition.
\item No two distinct signed branch words define the same affine map.
\item No two distinct signed branch words of length two define the same affine map.
\item Either $a_jb_j=0$ for every $j$, or $a_jb_{\tau(j)}=0$ for every $j$.
\end{enumerate}
\end{theorem}

\begin{proof}
A signed weak map at depth $k$ is determined by its geometric
$m$-adic interval and its final orientation. For two paths reading
the same geometric digits, the first distinct branch must take
opposite signs in the same original row. To finish in the same
orientation, the paths must later take opposite signs while their
orientations differ. At that step their original rows are a mirror
pair. Thus an exact overlap requires both a mixed row $h$ and
an index $j$ with $a_jb_{\tau(j)}>0$.
Conversely, these choices already give
\[
 g_{h,+}\circ g_{j,+}=g_{h,-}\circ g_{\tau(j),-}.
\]
This proves the equivalence of the last three conditions.

For each distinct signed map $g$ of depth $k$, let $w_k(g)$
count the planar words inducing it, and put $V_k=\sum_g w_k(g)^\beta$.
Grouping the signed maps by their geometric intervals gives
$Z_k\le V_k$, by $(x+y)^\beta\le x^\beta+y^\beta$.
Grouping pairs of signed maps by their composition gives, by the
same inequality, $V_{k+\ell}\le V_k V_\ell$. Since $V_1=U$, the
dimension formula gives the claimed upper bound. An exact overlap
at depth two makes the grouping at that depth strictly decrease the
sum, so $V_2<U^2$ and $\dimH K\le \log V_2/(2\log m)<\log U/\log m$.

In the absence of exact overlaps, every signed map has a unique
branch word, whose multiplicity is the product of the weights of
its letters. Consequently $V_k=U^k$. There are at most two signed
maps over each geometric interval. The inequality
$(x+y)^\beta\ge 2^{\beta-1}(x^\beta+y^\beta)$ thus gives
$Z_k\ge 2^{\beta-1}U^k$, proving equality in the dimension bound.

Finally, the open set condition precludes exact overlaps, since
distinct words of equal length take a feasible open set to disjoint
sets. Conversely, suppose there are no exact overlaps. If each row
has at most one signed branch then $(0,1)$ is feasible. Otherwise
choose a mixed row $h$. The last condition in the theorem forces
its mirror $\tau(h)$ to be empty. Let $I=(h/m,(h+1)/m)$ and set
\[
 O=\bigcup_{u\in\mathcal B^*}g_u(I),
\]
where the empty word is included. This open set is forward invariant.
To see disjointness, take words $u,v$ of positive length with distinct first letters,
and suppose $k=|u|\le |v|$. If $g_u(I)$ meets $g_v(I)$, then the
open level-$k$ cells of $u$ and the length-$k$ prefix $v'$ of $v$
coincide. Since $u\ne v'$ and there are no exact overlaps,
$g_{v'}=g_u\circ R$, where $R(y)=1-y$. Write $v=v'z$.
Intersection now requires $g_z(I)$ to meet $R(I)$, the open
interval in row $\tau(h)$. If $z$ is empty then $g_z(I)=I$, and
otherwise $g_z(I)$ lies in an active row. Both possibilities are
excluded because $\tau(h)$ is empty. Hence $O$ is feasible.
\end{proof}

\begin{corollary}\label{cor:block-osc}
Fix $k\ge1$. Identify identical weak maps at depth $k$ and retain one
copy of each. The resulting family satisfies the open set condition
if and only if the level-one signed branches do.
\end{corollary}

\begin{proof}
A feasible open set for the level-one family also works at every
depth. Conversely, if the level-one family fails the open set
condition, the theorem supplies a mixed row $h$ and a mirror pair
with $a_j b_{\tau(j)}>0$. Write $R(y)=1-y$.
Prepend any common word of length $k-1$ to $g_{h,+}$ and $g_{h,-}$.
The two resulting maps occupy the same level-$k$ interval and have
opposite signs, so this interval is a mixed row in base $m^k$.
Next append any common word of length $k-1$ to $g_{j,+}$ and
$g_{\tau(j),-}=R\circ g_{j,+}$. The resulting maps are $H$ and
$R\circ H$. They have opposite signs and occupy mirror rows in
base $m^k$. Both pairs survive identification of identical maps.
The length-two criterion of Theorem~\ref{thm:branches}, applied
in base $m^k$, rules out the open set condition at depth $k$.
\end{proof}

For example, if the occupied rows are disjoint from their mirrors, then $a_jb_{\tau(j)}=0$ for all $j$, and the signed-branch sum is exact. In particular, for a single occupied noncentral row with $a$ upright cells and $b$ upside-down cells the dimension is $\log_m(a^\beta+b^\beta)$. The classical value is $\log_n(a+b)$. Taking $a=b=1$ realizes the maximum gain in \eqref{eq:gain} on every grid.

\section{Box dimension and a spectral test for equality}\label{sec:box}

Two matrices of order two suffice to decide whether the Hausdorff and box
dimensions coincide. The first counts the possible weak-coordinate words.
The second counts pairs of paths reading the same word.

Put $d=\#D$, let $J=\{j:a_j+b_j>0\}$, and write
$R=\#J$, $B=\#\{j:a_jb_j>0\}$, $V=\#(J\cup\tau(J))$,
where $\tau(J)=\{\tau(j):j\in J\}$.
Allow both initial orientations. After reading $j$, the
numbers of paths ending in orientations $+$ and $-$ are
$a_j+b_{\tau(j)}$ and $b_j+a_{\tau(j)}$, respectively. Put
\[
 S=\#\{j:(a_j+b_{\tau(j)})(b_j+a_{\tau(j)})>0\},
 \qquad
 L=\begin{pmatrix}R-B&B\\V-S&S\end{pmatrix}.
\]
Define a second matrix, with all sums over $0\le j<m$, by
\[
 C=\begin{pmatrix}
 \sum_j(a_j^2+b_j^2)&2\sum_j a_jb_j\\
 \sum_j(a_jb_{\tau(j)}+b_ja_{\tau(j)})&
 \sum_j(a_ja_{\tau(j)}+b_jb_{\tau(j)})
 \end{pmatrix}.
\]

\begin{theorem}\label{thm:box}
Write $r=\rho(L)$ and $c=\rho(C)$ for the spectral radii. Then
\[
 \dim_{\rm B}K=\frac{\log d}{\log n}
       +\left(\frac1{\log m}-\frac1{\log n}\right)\log r,
 \qquad
 \dimH K=\dim_{\rm B}K\quad\Longleftrightarrow\quad rc=d^2.
\]
In all cases $rc\ge d^2$.
\end{theorem}

\begin{proof}
After a weak word has been read, the possible current orientations form
a singleton or the pair $\{+,-\}$. From a singleton, $R-B$ digits lead
to a singleton and $B$ to the pair. From the pair, the corresponding
counts are $V-S$ and $S$. Thus the number $u_k$ of geometric weak words
of length $k$ from one specified orientation is $(1,0)L^k(1,1)^t$.
When $B>0$ both types are accessible, so $u_k^{1/k}\to r$ by Perron--Frobenius theory. When $B=0$ each digit counted in $S$ must belong to $J\cap\tau(J)$, giving $S\le R$, so $u_k=R^k$ and $r=R$.
The two choices of initial orientation give mirror languages of the same size.

At scale $n^{-k}$ take $k$ digits in both coordinates and continue the
weak coordinate to depth $\ell=\lceil k\log n/\log m\rceil$.
Each of the $d^k$ distinct prefix rectangles admits exactly $u_{\ell-k}$
such continuations. The resulting grid cells have aspect ratios bounded
by $m$ and each contains a coded point. A ball of radius $n^{-k}$ meets
a bounded number of these cells, including in the presence of shared
boundaries. The covering number is therefore comparable to
$d^ku_{\ell-k}$, proving the box formula.

If every weak sign is positive, then $r=R$ and
$c=\sum_j a_j^2$ by Cauchy--Schwarz. Both the spectral equality and
equality of dimensions amount to constancy of the positive row counts,
by strict concavity in the classical row formula. Assume henceforth
that at least one weak sign is negative.

To identify the equality case, unfold both coordinate signs. The vertices
are the group generated by the sign pairs $(\sigma_e,\eta_e)$, $e\in D$,
under componentwise multiplication. At each vertex take one edge for each
$e$, multiply the orientation by its signs and label the edge by the
geometric grid digit at that orientation. This graph is strongly connected.
Its outgoing labels are distinct at each vertex. Its label shift $X$
therefore has at most one path per starting vertex for any given label word.
The geometric image of $X$ is the union of isometric copies of $K$
obtained by the group reflections. It has the same Hausdorff and box
dimensions as $K$, also upon passage to the torus, since the unit
square is covered by finitely many patches on which the quotient is
isometric. The full label shift and its weak factor $Y$ both satisfy
weak specification because two representing paths can be joined by a path of uniformly bounded length.

Choose the initial vertex uniformly and choose each successive edge
uniformly among its $d$ choices. The vertex law is stationary since
every generator permutes the vertices. Conditional on the initial vertex, the $d^k$ possible label words are distinct and equiprobable. Removing this conditioning changes the entropy by at most $\log 4$. The label process thus has entropy $\log d$, which is maximal since there are at most $4d^k$ words of length $k$. Weak specification gives uniqueness \cite[Proposition~2.2]{ZFeng}. Its weak projection $\mu$ satisfies
\[
 \mu([w])=d^{-k}vM_w\mathbf 1,\qquad
 v=(1/2,1/2),\quad |w|=k,
\]
because the weak sign is uniform in the group. The weak factor has
topological entropy $\log r$, since its word counts lie between $u_k$
and $2u_k$. The equality criterion for shifts with weak specification
\cite[Theorem~3.1]{ZFeng} now yields
\begin{equation}\label{eq:weak-max}
 \dimH K=\dim_{\rm B}K\quad\Longleftrightarrow\quad h(\mu)=\log r .
\end{equation}

Start two independent copies of this weak process. Conditional on their
current signs being equal or opposite, the entries of $C/d^2$ give the
probabilities of reading a common digit and ending with equal or opposite
signs. Hence
\[
 \sum_{|w|=k}\mu([w])^2
  =d^{-2k}(1/2,1/2) C^k\mathbf1.
\]
The initial and terminal vectors are strictly positive, so the $k$th root
of the right side tends to $c/d^2$, also for reducible or periodic $C$.
The Shannon entropy of a finite distribution dominates its order-two
R\'enyi entropy by Jensen's inequality. Passing to rates gives
$\log(d^2/c)\le h(\mu)\le\log r$. This proves $rc\ge d^2$, and
$rc=d^2$ implies equality of dimensions by \eqref{eq:weak-max}.

Conversely, if $h(\mu)=\log r$, the Gibbs bound for the maximal entropy
measure of a shift with weak specification gives
$\mu([w])\le A r^{-k}$ with a constant $A$ independent of $k,w$
\cite[Proposition~2.2]{ZFeng}. Therefore
$\sum_w\mu([w])^2\le A r^{-k}$, yielding $rc\le d^2$ as required.
\end{proof}

Feng, Lo and Shen gave a general test for dimension
coincidence for irreducible planar sofic affine-invariant sets~\cite{FLS}. The matrices
$L$ and $C$ above express the test directly in terms of the signed row counts.

For the classical carpets, Peres \cite{Peres} showed that unequal Hausdorff and box dimensions imply infinite Hausdorff measure in the critical dimension. For $s=\dimH K$,
the critical measure $\mathcal H^s(K)$ is positive and finite if $rc=d^2$.
If $rc>d^2$, then $\mathcal H^s$ on $K$ is not $\sigma$-finite.
Indeed, Qiu and Wang \cite[Theorem~1.1]{QW} give this dichotomy for the image of a shift with weak specification. The group graph above presents just such a shift for the union of the group images of $K$, also when every vertical sign is positive. For a single cell the assertion follows directly. Finite unions of isometric copies preserve both alternatives. Passage to the torus preserves them as well, by the same finite cover by isometric patches.

\section*{Acknowledgement}
The author thanks Mykola Pratsiovytyi for posing this question.
During the preparation of this work, the author used several large language models. The author reviewed and verified all AI-assisted content and takes full responsibility for the final manuscript.

\end{document}